\documentclass[11pt,twoside,a4paper]{article}
\usepackage{amsfonts}\usepackage{amsmath}\include{german}
\usepackage{amsthm}\usepackage{a4}
\usepackage{amssymb}\usepackage{graphics}\usepackage{url}
\numberwithin{equation}{section}
\newtheorem{thm}{Theorem}[section]
\newtheorem{cor}[thm]{Corollary}
\newtheorem{lem}[thm]{Lemma}
\theoremstyle{definition}

\newtheorem{rem}[thm]{Remark}
\newtheorem{exa}[thm]{Example}
\newcommand{\comment}[1]{}
\begin{document}
\title{Measure-preserving transformations of Volterra Gaussian processes and related bridges}\author{C\'eline Jost\footnote{Department of Mathematics and Statistics, P.O.\,Box 68 (Gustaf H\"allstr\"omin katu
2b), 00014 University of Helsinki, Finland.
 E-mail: celine.jost@helsinki.fi}}\date{}\maketitle
\begin{abstract}
We consider Volterra Gaussian processes on $[0,T]$, where $T>0$ is a fixed time horizon. These are processes of type $X_t=\int^t_0 z_X(t,s)dW_s$, $t\in[0,T]$,  where  $z_X$ is a square-integrable kernel and $W$ is a standard Brownian motion. An example is fractional Brownian motion. By using classical techniques from operator theory, we derive measure-preserving transformations of $X$ and their inherently related bridges of $X$. As a closely connected result, we obtain a Fourier-Laguerre series expansion for the first Wiener chaos of a Gaussian martingale over $[0,\infty)$.    
 \end{abstract}
\section{Introduction}
Fix $T>0$ and let $(X_t)_{t\in[0,T]}$ be a continuous \textit{Volterra} Gaussian process on a complete probability space $(\Omega,\mathcal{F},\mathbb{P})$. Hence, there exist a kernel $z_X \in L^2\left([0,T]^2\right)$ which is \textit{Volterra}, i.e. $z_X(t,s)=0,$ $s\geq t$, and a standard Brownian motion $(W_t)_{t\in[0,T]}$, such that 
\begin{equation}X_t\ =\ \int^t_0 z_X(t,s)dW_s,\ a.s.,\ t\in[0,T].\label{Volterra}\end{equation} Clearly, $X_0=0$, a.s., $X$ is centered and \begin{equation}R^X(s,t)\ =\ \int^{s\wedge t}_0 z_X(t,u)z_X(s,u)du,\ s,t\in[0,T],\label{covrep}\end{equation}
where $R^Y(s,t):=\text{Cov}_{\mathbb{P}}\left(Y_s,Y_t\right)$, $s,t\in[0,T]$, for a general Gaussian process $Y$.
We assume that $z_X$ is \textit{non-degenerate}, meaning that the family $\{z_X(t,\cdot)\,|\, t\in(0,T]\}$ is linearly independent and generates a dense subspace of $L^2([0,T])$.
On the one hand, it follows from this that $R^X$ is \textit{non-degenerate on} $(0,T]$, meaning that the matrices $\left\{ R^X(t_i,t_j)\right\}_{i,j=1,
\ldots,n}$, where $0<t_1<\ldots<t_n\leq T$ and $n\in\mathbb{N}$, are positive definite. On the other hand, it follows that \begin{equation}\Gamma_t(X)\ =\ \Gamma_t(W),\ t\in[0,T],\label{assumption0}\end{equation} i.e.
 \begin{equation}\mathbb{F}^{X}_{T}\ =\ \mathbb{F}^{W}_{T}.\label{filtrationX=filtrationW}\end{equation}
Here, $\Gamma_{T}(Y):=\overline{\text{span}\{Y_t\,|\,t\in[0,T]\}}\subseteq L^2(\mathbb{P})$ 
and $\mathbb{F}^{Y}_{T}:=\left(\mathcal{F}^Y_t\right)_{t\in[0,T]}$ denote the first Wiener chaos and the completed natural filtration of the Gaussian process $Y$ with $Y_0=0$, a.s., over $[0,T]$, respectively.
Moreover, we assume that  \begin{equation}z_X(T,\cdot)\ \neq\ 0,\ \lambda\text{-}a.e.,\label{kerinvertible}\end{equation}where $\lambda$ denotes Lebesgue measure, and
\begin{equation}z_X(\cdot,s)\ \text{has bounded variation on}\ [u,T]\ \text{for all}\ u\in(s,T).\label{kernelsingular}\end{equation} 
Processes of this type are a natural generalization of the nowadays in connection with finance and telecommunications extensively studied \textit{fractional Brownian motion with Hurst index} $H\in(0,1)$. Fractional Brownian motion, denoted by $B^H$, is the unique (in the sense of equality of finite-dimensional distributions) centered, $H$-self-similar Gaussian process with stationary increments. For $H=\frac{1}{2}$, fractional Brownian motion is standard Brownian motion. We have that \begin{equation*}z_{B^H}(t,s)\ =\ \frac{C(H)}{\Gamma\left(H+\frac{1}{2}\right)}(t-s)^{H-\frac{1}{2}} \cdot{}_2F_1\left(\frac{1}{2}-H,H-\frac{1}{2},H+\frac{1}{2},\frac{s-t}{s}\right)1_{[0,t)}(s),\label{kernelfBm}\end{equation*} 
$s,t\in[0,T]$, where $\Gamma$ denotes the gamma function, $C(H):=\left(\frac{2H\Gamma\left(H+\frac{1}{2}\right)\Gamma\left(\frac{3}{2}-H\right)}{\Gamma(2-2H)}\right)^{\frac{1}{2}}$ and ${}_2 F_1$ is the Gauss hypergeometric function. Moreover, \begin{equation*}R^{B^H}(s,t)\ =\ \frac{1}{2}\left(s^{2H}+t^{2H}-|t-s|^{2H}\right),\ s,t\in[0,T].\label{covariancefBm}\end{equation*}
For a Gaussian process $(Y_t)_{t\in[0,T]}$ with $Y_0=0$, a.s., recall that the process $\left(Y^T_t\right)_{t\in[0,T]}$ is a \textit{bridge of $Y$ (from $0$ to $0$ on $[0,T]$)} if  $$\text{Law}_{\mathbb{P}}\left(Y^T\right)\ =\ \text{Law}_{\mathbb{P}}\left(Y\, |\, Y_T=0\right).$$
 Clearly, $Y^T$ is a Gaussian process with $$Y^T_0\ =\ Y^T_T\ =\ 0,\ a.s.,$$ i.e. it is tied to $0$ at both ends. Therefore, the bridge is an intuitive concept for interpolation. It is well-known that the continuous process \begin{equation}\widehat{X}^T_t\ :=\ X_t\ -\ \frac{R^X(t,T)}{R^X(T,T)}X_T,\ t\in[0,T],\label{e:666}\end{equation} is a bridge of $X$ that satisfies \begin{equation}\Gamma_T\left(\widehat{X}^T\right)\ \perp\ \text{span}\{X_T\}\ =\ \Gamma_T(X),\label{notsamex}\end{equation} where $\perp$ denotes the orthogonal direct sum (see \cite{Ga}, Proposition 4).
$\widehat{X}^T$ is called the \textit{anticipative} bridge of $X$, since 
 \begin{equation}\mathbb{F}^{\widehat{X}^T}_T\vee\sigma(X_T)\ =\ \mathbb{F}^{X}_T\vee\sigma(X_T),\label{augmentedfiltrationscoincide}\end{equation} 
i.e. not the natural, but the initially $\sigma(X_T)$-enlarged filtrations of $X$ and $\widehat{X}^T$ coincide. Note that, more generally, $\widehat{Y}^T$ is a bridge of $Y$ for every Gaussian process $(Y_t)_{t\in[0,T]}$ with $Y_0=0$, a.s., and $R^Y(T,T)>0$.\\
\\
We implicitly assume that 
 $(\Omega,\mathcal{F},\mathbb{P})$ is the \textit{coordinate space of X}, which means that $\Omega=\{\omega:[0,T]\to\mathbb{R}\,|\,\omega\ \text{is continuous}\}$, $\mathcal{F}=\mathcal{F}^{X}_T$ and $\mathbb{P}$ is the probability measure with respect to which the \textit{coordinate process} $X_t(\omega)=\omega(t)$, $\omega\in\Omega$, $t\in [0,T]$, is a centered Gaussian process with covariance function $R^X$. $\mathcal{F}$ is also the Borel-$\sigma$-algebra on $\Omega$ equipped with the norm of uniform convergence.
The measurable map \begin{eqnarray*}\mathcal{T}\ :\ (\Omega,\mathcal{F},\mathbb{P})&\to&(\Omega,\mathcal{F},\mathbb{P})\\
X(\omega)\quad &\mapsto &\mathcal{T}(X(\omega))\end{eqnarray*}
 is a \textit{measure-preserving transformation on} $(\Omega,\mathcal{F},\mathbb{P})$ if it is a endomorphism, i.e. if $\mathbb{P}^{\mathcal{T}}=\mathbb{P}$, or equivalently, if  $\mathcal{T}(X)\stackrel{d}{=}X$, where $\stackrel{d}{=}$ denotes equality of finite-dimensional distributions.
The measurable map \begin{eqnarray*}\mathcal{B}\ :\ (\Omega,\mathcal{F},\mathbb{P})&\to&(\Omega,\mathcal{F},\mathbb{P})\\
X(\omega)\quad &\mapsto &\mathcal{B}(X(\omega))\end{eqnarray*}
 is a \textit{bridge transformation on} $(\Omega,\mathcal{F},\mathbb{P})$ if $\mathcal{B}(X)\stackrel{d}{=}X^T$.
Note that, in particular, $\mathcal{T},\mathcal{B}:C([0,T])\to C([0,T])$ are linear maps. \\
\\
In this work, we derive measure-preserving transformations, such that\begin{equation*}\Gamma_T(\mathcal{T}(X))\ =\ \Gamma_{T}\left(\widehat{X}^T\right).\label{introTX=ANTX}\end{equation*}
As a naturally related problem, we derive bridge transformations that satisfy
\begin{equation*}\Gamma_T(\mathcal{B}(X))\ =\ \Gamma_T(X).\label{introBX=X?}\end{equation*}  
We define $\mathcal{T}(X)$ and $\mathcal{B}(X)$ as Wiener integral processes with respect to $X$, where the integrand kernels are based on suitable isometric operators on the Wiener integrand space of $X$. This method follows ideas of Jeulin and Yor (see \cite{Je}) and Peccati (see \cite{Pe}), where the case $X=W$ is considered. If $X$ is a martingale, then due to independence of increments, the operators and the resulting transformations are simple. However, if $X$ is not a martingale, then operators and obtained transformations are technically more involved. In fact, the operators are composed by the corresponding operators for the prediction martingale of $X_T$, and the isometry between the Wiener integrand spaces of $X$ and the prediction martingale, respectively. From these results, we obtain in particular an alternative, purely operator theoretic derivation for the \textit{dynamic} bridge of $X$, as introduced in \cite{Ga}.
Furthermore, by using the same operator theoretic methods and as a closely connected result, we obtain a
Fourier-Laguerre series expansion for the first Wiener chaos of a Gaussian martingale over $[0,\infty)$.
\begin{rem}\label{canonicalrepresentation}1. Representation (\ref{Volterra}) is unique in the following sense: assume that there exist another non-degenerate Volterra kernel $z'_X$ and another standard Brownian motion $W'$, such that $X_t=\int^t_0z_X'(t,s)dW'_s$, a.s., $t\in[0,T]$. Then, it follows from (\ref{assumption0})  that $\Gamma_t(W)=\Gamma_t(W')$, $t\in[0,T]$, i.e. $W$ and $W'$ are indistinguishable. In particular, $$0\ =\ \text{E}_{\mathbb{P}}(X_t-X_t)^2\ =\ \int^t_0 \left(z_X(t,s)-z_X'(t,s)\right)^2 ds,\ t\in[0,T],$$ which implies that $z_X(t,\cdot)=z_X'(t,\cdot)$, $\lambda$\text{-}a.e., $t\in[0,T]$.\\
2. Although the covariance function $R^X$ is continuous on $[0,T]^2$, it is  generally not true that $z_X(\cdot,s)$ is continuous for all $s\in[0,T]$. For example, for the fractional Brownian motion with $H<\frac{1}{2}$, the function $z_{B^H}(\cdot,s)$ is discontinuous in $t=s$ for every $s\in(0,T)$.\\
3. (\ref{kerinvertible}) is not a necessary condition for the non-degeneracy of a Volterra kernel as the example $z_X(t,s):=1_{[\frac{t}{2},t)}(s)$, $s,t\in[0,T]$, shows.
\end{rem}
The article is organized as follows. Section 2 is dedicated to Wiener integration with respect to $X$ and $X^T$, respectively. In Section 3, we derive the measure-preserving and  bridge transformations in the case, when $X$ is a martingale. These results are generalized to Volterra Gaussian processes in Section 4. In Section 5, we derive the series expansion.
\section{Wiener integrals} First, we review the construction of abstract and time domain Wiener integrals with respect to the Volterra Gaussian process $X$. Second, we define abstract Wiener integrals with respect to the bridge $X^T$, and explain their relation to Wiener integrals with respect to $X$. \subsection{Wiener integrals with respect to $X$}
Let $\mathcal{E}_T:=\text{span}\{1_{[0,t)}\, |\,t\in(0,T]\}$ be the space of elementary functions on $[0,T]$.  The Hilbert space of \textit{abstract Wiener integrands of $X$ on $[0,T]$}, denoted by $\Lambda_T(X)$, is defined as  the completion of $\mathcal{E}_T$ with respect to the scalar product $$(1_{[0,s)},1_{[0,t)})_{X}\ :=\ R^X(s,t),\ s,t\in(0,T].$$ %
Hence, $f\in\Lambda_T(X)$ is an equivalence class of Cauchy sequences $\{f_n\}_{n\in\mathbb{N}}\in\mathcal{E}_T^{\mathbb{N}}$, where $\{f_n\}_{n\in\mathbb{N}}\sim \{g_n\}_{n\in\mathbb{N}}:\Leftrightarrow (f_n-g_n,f_n-g_n)_X \to 0$, $n\to\infty$. The scalar product  on $\Lambda_T(X)$ is given by $$(f,g)_{\Lambda_T(X)}:=\lim_{n\to\infty}(f_n,g_n)_X,\ f,g\in\Lambda_T(X),\  \{f_n\}_{n\in\mathbb{N}}\in f,\ \{g_n\}_{n\in\mathbb{N}}\in g,$$ and induces the norm $|\cdot|_{\Lambda_T(X)}$. 
The  isometric isomorphism defined by  \begin{eqnarray*}I^X_T\ :\ \Lambda_T(X)&\to&\Gamma_T(X)\\
1_{[0,t)}\ \ &\mapsto &X_t,\ t\in(0,T],\end{eqnarray*} is  called the \textit{abstract Wiener integral with respect to $X$ on $[0,T]$}.  By definition, $\int^T_0 f(s)dX_s:=I^X_T(f)$ is centered, Gaussian and $\text{E}_{\mathbb{P}}\left(\int^T_0 f(s)dX_s\right)^2=|f|^2_{\Lambda_T(X)}$ for all $f\in\Lambda_T(X)$.\\
\\
By combining  (\ref{Volterra}), which is also called the \textit{time domain representation of} $X$, with the standard Wiener integral, we obtain a subspace of $\Lambda_T(X)$ whose elements can be identified with functions. For this purpose, define a linear operator $$\left(\mathbf{K}^X f\right)(s)\ :=\ f(s)z_X(T,s)\ +\ \int^T_s \left(f(u)-f(s)\right)z_X(du,s),\ s\in(0,T).$$
By (\ref{kernelsingular}), the integral is well-defined for a suitably large class of functions. Clearly, $\mathbf{K}^X$ extends the linear isometry \begin{eqnarray*}\left(\mathcal{E}_T,\left(\cdot\, ,\cdot\right)_X\right)\ &\to&L^2([0,T])\\ 
1_{[0,t)}\qquad&\mapsto&z_{X}(t,\cdot),\ t\in(0,T].\end{eqnarray*} We define the space of \textit{time domain Wiener integrands of $X$ on $[0,T]$} by  \begin{equation*}\check{\Lambda}_T(X)\ :=\ \left\{ f:[0,T]\to\mathbb{R}\ \bigg|\ \mathbf{K}^{X}f\ \text{is well-defined and}\ \int^T_0\left(\mathbf{K}^{X}f\right)^2(s)ds<\infty\right\}\end{equation*} with  scalar product $$\left(f,g \right)_{\check{\Lambda}_T(X)}\ :=\ \left(\mathbf{K}^X f ,\mathbf{K}^X g\right)_{L^2([0,T])},\ f,g\in\check{\Lambda}_T(X).$$ 
By (\ref{covrep}), we have that $$(1_{[0,s)},1_{[0,t)})_{\check{\Lambda}_T(X)}\ =\ (1_{[0,s)},1_{[0,t)})_{\Lambda_T(X)},\ s,t\in[0,T].$$
Let $f\in\Lambda_T(X)$. If there exists $\check{f}\in\check{\Lambda}_T(X)$, such that 
$$\left(f,1_{[0,t)}\right)_{\Lambda_T(X)}\ =\ \left(\check{f},1_{[0,t)}\right)_{\check{\Lambda}_T(X)}\ \text{for all}\ t\in(0,T],$$ then we identify $f$ and $\check{f}$.
$\mathcal{E}_T$ is dense in $\check{\Lambda}_T(X)$, hence $\check{\Lambda}_T(X)\subseteq \Lambda_T(X)$. However, in general, $\check{\Lambda}_T(X)\not = \Lambda_T(X)$.  The restriction $I^X_T|_{\check{\Lambda}_T(X)}$ is called the \textit{time domain Wiener integral with respect to $X$ on $[0,T]$}.
\begin{rem}\label{normf0ifff0}Due to the non-degeneracy of $z_X$, we have that $|f|_{\check{\Lambda}_T(X)}=0$ if and only if $f=0$, $\lambda$-a.e.\end{rem}
\begin{rem}From (\ref{assumption0}), it follows that (\ref{Volterra}) has a reciprocal in the sense that there exists a Volterra kernel $z_{X}^{\ast}$ with $z^{\ast}_{X}(t,\cdot)\in\Lambda_T(X)$, $t\in(0,T]$, such that \begin{equation*}W_t\ =\ \int^t_0 z^{\ast}_X(t,s)dX_s,\ a.s.,\ t\in[0,T].\end{equation*}
For $t\in(0,T]$, this abstract Wiener integral is a time domain Wiener integral if and only if there exists $f\in\check{\Lambda}_{T}(X)$, such that $\mathbf{K}^X f =1_{[0,t)}$, $\lambda$-a.e.\end{rem}
\begin{rem}
If $R^X$ is of bounded variation, then it determines a finite signed measure on $[0,T]^2$. Hence, one can define an alternative space of Wiener integrands of $X$ on $[0,T]$ by $$|\Lambda_T(X)|\ :=\ \left\{f:[0,T]\to\mathbb{R}\ \bigg|\ \int^T_0\int^T_0 |f(s)||f(t)|\big|R^X\big|(dt,ds)<\infty\right\},$$ where $\big|R^X\big|$ denotes the measure of total variation of $R^X$, and the integral is a Lebesgue-Stieltjes integral. The corresponding scalar product is given by $$\left( f,g\right )_{|\Lambda_T(X)|}\ :=\ \int^T_0 \int^T_0 f(t)g(s)R^{X}(dt,ds),\ f,g\in|\Lambda_T(X)|.$$ For details, see \cite{Hu}.  
\end{rem}
\begin{exa}For the fractional Brownian motion holds that (see \cite{Al}, p. 797-800)
 $$\left(\mathbf{K}^{B^H}f\right)(s)\ =\ C(H)s^{\frac{1}{2}-H}\left(\mathcal{I}_{T-}^{H-\frac{1}{2}}\cdot^{H-\frac{1}{2}}f\right)(s),\ s\in(0,T),$$ where $\mathcal{I}_{T-}^{H-\frac{1}{2}}$ denotes the right-sided Riemann-Liouville fractional integral operator of order $H-\frac{1}{2}$ over $[0,T]$. Furthermore, 
$$z^{\ast}_{B^H}(t,s)\ =\
\frac{C(H)^{-1}}{\Gamma\left(\frac{3}{2}-H\right)}(t-s)^{\frac{1}{2}-H}\cdot {}_2 F_1\left(\frac{1}{2}-H,\frac{1}{2}-H,\frac{3}{2}-H,\frac{s-t}{s}\right)1_{[0,t)}(s),$$ $s,t\in[0,T]$.
 We have that
 $\check{\Lambda}_T\left(B^H\right)=\Lambda_T\left(B^H\right)$ if and only if $H\leq\frac{1}{2}$.   
 Clearly, $R^{B^H}$ is of bounded variation if and only if $H\geq\frac{1}{2}$. Then,
  $z_{B^H}^{\ast}(t,\cdot)\in\big|\Lambda_T\left(B^H\right)\big|\subseteq\check{\Lambda}_T\left(B^H\right)$, $t\in(0,T]$.
 See \cite{Pi} for full proofs.\end{exa}
\subsection{Wiener integrals with respect to $X^T$}
From (\ref{e:666}), it follows that $X^T$ is centered and \begin{equation}R^{X^T}\!(s,t) \ =\ R^X(s,t)\ -\ \frac{R^X(s,T)R^X(t,T)}{R^X(T,T)},\ s,t\in[0,T].\label{covbridge}\end{equation}
Hence,  $R^{X^T}$ is non-degenerate on $(0,T)$ and $R^{X^T}\!(T,T)=0$. In order to define a scalar product based on $R^{X^T}$, let $$\mathcal{C}_T\ :=\ \{f:[0,T]\to\mathbb{R}\ |\ f\ \text{is}\ \text{constant}\ \lambda\text{-}a.e.\}$$ and consider the quotient space $\mathcal{E}_T/\mathcal{C}_T$. Clearly,  $1_{[0,T)}\sim 0$ mod $\mathcal{C}_T$.
 The Hilbert space of \textit{abstract Wiener integrands of  $X^T$ on $[0,T]$}, denoted by $\Lambda_T\!\left(X^T\right)$, is defined as the completion of $\mathcal{E}_T/\mathcal{C}_T$ with respect to the scalar product 
$$\left(\overline{1_{[0,s)}},\overline{1_{[0,t)}}\right)_{X^T} \  :=\ R^{X^T}\!(s,t),\ s,t\in(0,T].$$
Thus, $\overline{f}\in\Lambda_T\!\left(X^T\right)$ is an equivalence class of Cauchy sequences $\left\{\overline{f_n}\right\}_{n\in\mathbb{N}}\in(\mathcal{E}_T/\mathcal{C}_T)^{\mathbb{N}}$, where $\left\{\overline{f_n}\right\}_{n\in\mathbb{N}}\sim \{\overline{g_n}\}_{n\in\mathbb{N}}:\Leftrightarrow \left(\overline{f_n-g_n},\overline{f_n-g_n}\right)_{X^T} \to 0$, $n\to\infty$. The scalar product on $\Lambda_T\!\left(X^T\right)$ is given by $$\left(\overline{f},\overline{g}\right)_{\Lambda_T(X^T)}:=\lim_{n\to\infty}\left(\overline{f_n},\overline{g_n}\right)_{X^T},\  \overline{f},\overline{g}\in\Lambda_T\!\left(X^T\right),\ \left\lbrace\overline{f_n}\right\rbrace_{n\in\mathbb{N}}\in \overline{f},\ \{\overline{g_n}\}_{n\in\mathbb{N}}\in\overline{g},$$ with induced norm $|\cdot|_{\Lambda_T\left(X^T\right)}$.
The  isometric isomorphism defined by \begin{eqnarray*}I^{X^T}_T\ :\ \Lambda_T\!\left(X^T\right)&\to&\Gamma_T\!\left(X^T\right)\\
\overline{1_{[0,t)}}\quad &\mapsto &X^T_t,\ t\in(0,T],\end{eqnarray*} is  called the \textit{abstract Wiener integral with respect to $X^T$ on $[0,T]$.}
We have that $\int^T_0 \overline{f}(s)dX^T_s$ $:=I^{X^T}_T\!\left(\overline{f}\right)$ is centered, Gaussian and  $\text{E}_{\mathbb{P}}\left(\int^T_0 \overline{f}(s)dX^T_s\right)=\big|\overline{f}\big|^2_{\Lambda_T\left(X^T\right)}$ for all $\overline{ f}\in\Lambda_T\!\left(X^T\right)$.\\
\\
 Let \begin{eqnarray*}\eta^X\ :\ \Lambda_T(X) &\to&\Lambda_T(X)\nonumber\\
f(\cdot)\quad &\mapsto&  f(\cdot)\ -\ \frac{\left(f,1_{[0,T)}\right)_{\Lambda_T(X)}}{R^X(T,T)}.\label{etta}\end{eqnarray*} 
$\eta^X$ is the linear orthoprojection from $\Lambda_T(X)$ onto the closed subspace $$\Lambda_{T,0}(X)\ :=\ \left\lbrace f\in \Lambda_T(X)\ \bigg|\ \left(f,1_{[0,T)}\right)_{\Lambda_T(X)}=0\right\rbrace$$ along $\Lambda_{T,0}(X)^{\perp}=\mathcal{C}_T$. In particular, $\text{ker}\left(\eta^X\right)=\mathcal{C}_T$ and \begin{equation}\eta^X(\Lambda_T(X))\ =\ \Lambda_{T,0}(X).\label{etalambda=lambda0}\end{equation}  Thus, the map\begin{eqnarray*}\overline{\eta^X}\ :\ \Lambda_T(X)/\mathcal{C}_T &\to& \Lambda_{T,0}(X) \nonumber\\
\overline{f}\qquad &\mapsto &\eta^X f\label{isomorphism}\end{eqnarray*}
is an isomorphism. 
From (\ref{covbridge}), we obtain that
\begin{equation*}\left(\overline{1_{[0,s)}},\overline{1_{[0,t)}}\right)_{X^T}\ =\ \left(\eta^X 1_{[0,s)},\eta^X 1_{[0,t)}\right)_{X},\ s,t\in(0,T].\label{ecip}\end{equation*}
Hence, \begin{equation*} \Lambda_T\!\left(X^T\right)\ =\ \Lambda_T(X)/\mathcal{C}_T,\label{spaces}\end{equation*} i.e. $$\overline{f}\in\Lambda_T\left(X^T\right)\ \Leftrightarrow\ f\in\Lambda_T(X).$$
From (\ref{e:666}), it follows that \begin{eqnarray}\widehat{X}^T_t &=& \int^T_0 \left(\eta^X 1_{[0,t)}\right)(s)dX_s,\ a.s., \ t\in[0,T].\label{pabr}\end{eqnarray}
Hence,
\begin{equation}\int^T_0 \overline{f}(s)d\widehat{X}^{T}_s\ =\ \int^T_0 \left(\eta^X f\right)(s)dX_s, \ a.s.,\  f\in\Lambda_T(X).\label{notinvx}\end{equation}
From (\ref{notsamex}), we see that (\ref{pabr}) does not have a reciprocal, i.e. $X_t$ can not be written as an abstract Wiener integral with respect to $\widehat{X}^T$. By using (\ref{notinvx}), we can write (\ref{notsamex}) as follows: \begin{equation}\int^T_0 f(s)dX_s\  =\ \int^T_0 \overline{f}(s)d\widehat{X}^T_s\  +\ \frac{\left(f,1_{[0,T)}\right)_{\Lambda_T(X)}}{R^X(T,T)}\cdot  X_T,\ f\in \Lambda_T(X).\label{intX=intantX+R}\end{equation}
By combining (\ref{notinvx})  and (\ref{etalambda=lambda0}),  we obtain the following:
\begin{lem}\label{chaosofantbridge=}We have that
 $$\Gamma_T\!\left(\widehat{X}^T\right)\ =\ I^{X}_T\left(\Lambda_{T,0}(X)\right).$$   \end{lem}
\section{Gaussian martingales}\label{section2}We consider the special case when $X$ is an $\mathbb{F}^{X}_T$-martingale (or equivalently, an $\mathbb{F}^W_T$-martingale). This corresponds to a Volterra kernel of type $z_X(t,s)=z_M(s)1_{[0,t)}(s)$, $s,t\in[0,T]$, where $z_M\in L^2([0,T])$ and $z_M\neq 0$, $\lambda$-a.e.  For convenience, we write $X=M$.
Let $\langle M\rangle_{\cdot}:=\int^{\cdot}_0 z^2_M(s)ds$ denote the quadratic variation or variance function of $M$. Then $\langle M\rangle_{\cdot}$ is strictly increasing on $[0,T]$. We denote \begin{equation*}\langle M\rangle_{T,t}\ :=\ \langle M\rangle_T\ -\ \langle M\rangle_t,\ t\in[0,T].\label{vardiff}\end{equation*}
We derive two types of measure-preserving and bridge transformations, present some pathwise relations between the transformed processes, and show how the two types a connected. 
For the case $M=W$, large parts of the results have been obtained in \cite{Je} and \cite{Pe}, although not in the present form and not so detailed. The generalization from $W$ to $M$ can be obtained by using the fact that $M\stackrel{d}{=}W_{\langle M\rangle_{\cdot}}$. However, for convenience, we  provide independent and complete proofs. 
\subsection{Measure-preserving and bridge transformations}
Clearly, we have that
 $$\Lambda_T(M)\ =\ L^2_T(M)\ :=\ L^2([0,T],d\langle M\rangle_{\cdot}).$$  Furthermore,  $$\Lambda_{T,0}(M)\ =\ L^2_{T,0}(M)\ :=\ \left\{f\in L^2_T(M)\ \bigg|\ \int^T_0 f(s)d\langle M\rangle_s=0\right\}$$ and
 $$\left(\eta^M f\right)(\cdot)\ = \ f(\cdot)\ -\ \frac{1}{\langle M\rangle_T}\int^T_0 f(s)d\langle M\rangle_s,\ f\in L^2_{T}(M).$$
\begin{lem}\label{pec}1. For $f\in L^2_T(M)$, let 
$$\left(\mathcal{H}^{M,1}f\right)(s)\ := \frac{1}{\langle M\rangle_s}\int^s_0 f(u)d\langle M\rangle_u,\ s\in(0,T],$$ and $$\left(\mathcal{H}^{M,2} f\right)(s)\ :=\    \frac{1}{\langle M\rangle_{T,s}}\int^T_{s} f(u)d\langle M\rangle_u,\ s\in[0,T).$$ The Hardy type operators $\mathcal{H}^{M,1}$ and $\mathcal{H}^{M,2}$ are bounded endomorphisms on $L^2_T(M)$ with adjoints $$\left(\mathcal{H}^{M,1,\ast}f\right)(s)\ =\  \int^T_{s} \frac{f(u)}{\langle M\rangle_u}d\langle M\rangle_u,\ s\in(0,T],$$
and
$$ \left(\mathcal{H}^{M,2,\ast}f\right)(s)\ =\ \int^{s}_0\frac{ f(u)}{\langle M\rangle_{T,u}}d\langle M\rangle_u, \ s\in[0,T),$$ respectively.\\
2. For $i\in\{1,2\}$, let furthermore  $$\alpha^{M,i} f\ :=\ f\ -\ \mathcal{H}^{M,i} f$$ and
  $$\beta^{M,i} f\ :=\  f\ -\ \mathcal{H}^{M,i,\ast}f.$$  Then
$\alpha^{M,i}:L^2_{T,0}(M)\to L^2_T(M)$ and $\beta^{M,i}:L^2_T(M)\to L^2_{T,0}(M)$ are isometric 
 isomorphisms with $\left(\alpha^{M,i}\right)^{-1}=\beta^{M,i}$. \\
3. For $i\in\{1,2\}$, we have that \begin{equation}\alpha^{M,i}\eta^{M} f\ =\ \alpha^{M,i} f,\ f\in\ L^2_T(M),\label{ea=a}\end{equation} and \begin{equation}\eta^{M}\beta^{M,i} f\ =\ \beta^{M,i} f,\ f\in\ L^2_T(M).\label{eb=b}\end{equation} Also,  \begin{equation}\beta^{M,i}\alpha^{M,i} f\ = \ \eta^M f, \ f\in L^2_T(M).\label{ba=e}\end{equation}\end{lem}
\begin{proof} 1. From Hardy's inequality (see \cite{Ha}, Theorem 327, p. 240), it follows that
 \begin{equation}\int^1_0\frac{1}{x^2}\left(\int^x_0 g(z)dz\right)^2dx \leq\ 4\int^1_0 g^2(x)dx, \ g\in L^2([0,1]).\label{hardy2}\end{equation}
By using  (\ref{hardy2}) with $g(z):=f\left(y^{-1}(z)\right)$ for $i=1$, and $g(z):=f\left(y^{-1}(1-z)\right)$ for $i=2$, where $y(u):=\frac{\langle M \rangle_u}{\langle M\rangle_T}$, we obtain that $\big|\mathcal{H}^{M,i}f\big|^2_{L^2_T(M)}\leq  4|f|^2_{L^2_T(M)}$, $f\in L^2_T(M)$.    
From Fubini's theorem, it follows that $\mathcal{H}^{M,i,\ast}$ is  the adjoint of $\mathcal{H}^{M,i}$. \\
2. By using Fubini's theorem and splitting integrals, we have that $\alpha^{M,i}:L^2_{T,0}(M)\to L^2_T(M)$ is the inverse of $\beta^{M,i}: L^2_T(M) \to L^2_{T,0}(M)$. Moreover, from part 1, it follows that $\alpha^{M,i}$ is the adjoint of $\beta^{M,i}$. So $\alpha^{M,i}$ and $\beta^{M,i}$ are unitary  and hence isometric.\\
3. A straightforward calculation yields (\ref{ea=a}). (\ref{eb=b}) follows from part 2. By using Fubini's theorem and splitting integrals, we obtain (\ref{ba=e}).\end{proof}
\begin{rem}\label{oprem}
1. Let $f\in L^2_T(M)$. For $t\in(0,T]$ and $s\in[t,T)$, we have that $\left(\alpha^{M,2}f 1_{[0,t)}\right)(s)=0$ and  $\left(\beta^{M,2}f 1_{[0,t)}\right)(s)=c_M(f,t)$, where $c_M$ is a function independent of $s$.\\
2. It holds that $\alpha^{M,i}c\equiv 0$, $c\in\mathcal{C}_T$, $i\in\{1,2\}$.\end{rem}
We can state and prove the following theorem:
\begin{thm}\label{tablem}Let $i\in\{1,2\}$. The transformation defined by \begin{equation}\mathcal{T}^{(i)}_t(M)\ :=\ \int^T_0\left(\beta^{M,i} 1_{[0,t)}\right)(s)dM_s,\ t\in [0,T],\label{pwc2}\end{equation} 
 is measure-preserving, i.e. $\mathcal{T}^{(i)}(M)$ is an $\mathbb{F}_{T}^{\mathcal{T}^{(i)}(M)}$-martingale with $\mathcal{T}^{(i)}_0(M)=0$, a.s., and $\langle\mathcal{T}^{(i)}(M)\rangle_{\cdot}=\langle M\rangle_{\cdot}$. Furthermore, the
process
 \begin{equation}\mathcal{B}^{(i)}_t(M)\  :=\ \int^T_0 \left(\alpha^{M,i}1_{[0,t)}\right)(s)dM_s,\ t\in[0,T],\label{pwc1}\end{equation} is a bridge of $M$. 
 \end{thm}
\begin{proof}Clearly, $\mathcal{T}^{(i)}(M)$ and $\mathcal{B}^{(i)}(M)$ are centered Gaussian processes. From part 2 of Lemma \ref{pec}, we obtain that $$\text{Cov}_{\mathbb{P}}\left(\mathcal{T}^{(i)}_s(M),\mathcal{T}^{(i)}_t(M)\right)\ =\ \text{Cov}_{\mathbb{P}}(M_s,M_t),\ s,t\in[0,T].$$
Furthermore, by combining  (\ref{ea=a}), part 2 of Lemma \ref{pec} and (\ref{pabr}), we obtain that  $$\text{Cov}_{\mathbb{P}}\left(\mathcal{B}^{(i)}_s(M),\mathcal{B}^{(i)}_t(M)\right)\ =\ \text{Cov}_{\mathbb{P}}\left(M^T_s,M^T_t\right),\ s,t\in[0,T].\qedhere$$ \end{proof}
\subsection{Pathwise relations for the transformed processes}\label{subsection3.2}
By using the stochastic Fubini theorem, we have that \begin{equation}\mathcal{T}^{(1)}_t(M)\ =\ M_t\ -\ \int^t_0 \frac{M_s}{\langle M\rangle_s}d\langle M\rangle_s,\ a.s.,\ t\in[0,T].\label{genyor}\end{equation}  Similarly, we obtain that \begin{equation}\mathcal{T}^{(2)}_t(M)\ =\ M_t\ -\ \int^t_0\frac{M_T-M_s}{\langle M\rangle_{T,s}}d\langle M\rangle_s,\ a.s.,\ t\in[0,T].\label{T^2(M)=}\end{equation}
On the one hand, by combining (\ref{pwc2}), (\ref{eb=b}) and (\ref{notinvx}), we have that \begin{eqnarray}\mathcal{T}^{(i)}_t(M)
&=&\int^T_0 \left(\overline{\beta^{M,i} 1_{[0,t)}}\right)(s)d\widehat{M}^T_s,\ a.s., \ t\in[0,T],\ i\in\{1,2\}.\label{T=intAB}\end{eqnarray}
On the other hand, by combining (\ref{pabr}), (\ref{ba=e}) and (\ref{pwc2}), we obtain the  reciprocal \begin{eqnarray}\widehat{M}^T_t 
&=& \int^T_0 \left(\alpha^{M,i}1_{[0,t)}\right)(s)d\mathcal{T}^{(i)}_s(M),\ a.s.,\ t\in[0,T],\ i\in\{1,2\}.\label{AB=intT}\end{eqnarray}
Clearly, from (\ref{T=intAB}) and (\ref{AB=intT}), it follows that \begin{equation}\Gamma_T\!\left(\mathcal{T}^{(i)}(M)\right)\ =\ \Gamma_T\!\left(\widehat{M}^T\right),\ i\in\{1,2\}.\label{GT=GANTB}\end{equation}
For $i=2$, it follows more precisely from part 1 of Remark \ref{oprem} that \begin{equation}\Gamma_t\!\left(\mathcal{T}^{(2)}(M)\right)\ =\ \Gamma_t\!\left(\widehat{M}^T\right),\ t\in[0,T],\label{gammatT2M=gammatANTM}\end{equation} i.e. \begin{equation}\mathbb{F}^{\mathcal{T}^{(2)}(M)}_{T}\ =\ \mathbb{F}^{\widehat{M}^T}_{T}.\label{F^T2(M)=}\end{equation}
From (\ref{F^T2(M)=}), we obtain that $\mathcal{T}^{(2)}(M)$ is an $\mathbb{F}^{\widehat{M}^T}_T$-martingale. Furthermore, by combi\-ning (\ref{GT=GANTB}) and (\ref{notsamex}), we have that $\mathcal{T}^{(2)}(M)$ is orthogonal to $M_T$, and thus independent of $M_T$. Hence, $\mathcal{T}^{(2)}(M)$ is also an $\left(\mathbb{F}^{\widehat{M}^T}_T\vee\sigma(M_T)\right)$-martingale and so, by 
using (\ref{augmentedfiltrationscoincide}), an  $\left(\mathbb{F}^{M}_T\vee\sigma(M_T)\right)$-martingale. From (\ref{T^2(M)=}), it is then easy to see that $\mathcal{T}^{(2)}(M)$ is the martingale component in the $\left(\mathbb{F}_T^{M}\vee \sigma(M_T)\right)$-semimartingale decomposition of $M$.\\
\\
By combining (\ref{T^2(M)=}) and (\ref{e:666}), we obtain that   
\begin{eqnarray*}\mathcal{T}^{(2)}_t(M)
&=&\widehat{M}^T_t\ +\ \int^t_0 \frac{\widehat{M}^T_s}{\langle M\rangle_{T,s}}d\langle M\rangle_s,\ a.s.,\ t\in[0,T].\end{eqnarray*}
In particular, $\widehat{M}^T$ is the unique solution of the linear stochastic differential equation $$d\widehat{M}^T_t\ =\ d\mathcal{T}^{(2)}_t(M)\ -\ \frac{\widehat{M}^T_t}{\langle M\rangle_{T,t}}d\langle M\rangle_t,\qquad \widehat{M}^{T}_0\ =\ 0.$$
A straightforward calculation yields \begin{equation}\mathcal{B}^{(1)}_t(M)\ =\ -\langle M\rangle_t\int^T_t \frac{1}{\langle M\rangle_s}dM_s,\ a.s.,\ t\in(0,T].\label{B1M=}\end{equation} Similarly, \begin{equation*}\mathcal{B}^{(2)}_t(M)\ =\ \langle M\rangle_{T,t}\int^t_0\frac{dM_s}{\langle M\rangle_{T,s}},\ a.s.,\ t\in[0,T).\label{B2M=}\end{equation*}
By combining part 2 of Lemma \ref{pec} and (\ref{pwc1}), we obtain that \begin{equation}M_t\ =\ \int^T_0 \left(\overline{\beta^{M,i}1_{[0,t)}}\right)(s)d\mathcal{B}^{(i)}_s(M), \ a.s.,\ t\in[0,T], \ i\in\{1,2\}.\label{betaintegral}\end{equation}
By comparing identities (\ref{pwc1}) and (\ref{betaintegral}) with identities (\ref{AB=intT}) and (\ref{T=intAB}), we observe that the bridge $\mathcal{B}^{(i)}(M)$ is related to the martingale $M$ in the same way as
 the bridge $\widehat{M}^T$ is related to the martingale $\mathcal{T}^{(i)}(M)$. Hence, from (\ref{GT=GANTB}), we obtain that \begin{equation}\Gamma_T\!\left(\mathcal{B}^{(i)}(M)\right)\ =\ \Gamma_T(M),\ i\in\{1,2\}.\label{gammabm=gammam}\end{equation}Moreover, from (\ref{gammatT2M=gammatANTM}) we have that $$\Gamma_t\!\left(\mathcal{B}^{(2)}(M)\right)\ =\ \Gamma_t(M),\  t\in[0,T],$$ i.e. $$\mathbb{F}^{\mathcal{B}^{(2)}(M)}_{T}\ =\ \mathbb{F}^{M}_{T}.$$
From (\ref{AB=intT}) and (\ref{betaintegral}), we obtain the following: 
\begin{lem}\label{lemt}Let $i\in\{1,2\}$. Then
 \begin{equation*}  \widehat{M}_t^{T} =\ \mathcal{B}^{(i)}_t\left(\mathcal{T}^{(i)}(M)\right),   \ a.s., \ t\in[0,T].\end{equation*}
Furthermore,  $$M_t\ =\ \mathcal{T}^{(i)}_t\left(\mathcal{B}^{(i)}(M)\right),\ a.s.,\ t\in[0,T].$$\end{lem}
\begin{rem}1. From (\ref{GT=GANTB}), we have that $\mathcal{F}^{\mathcal{T}^{(i)}(M)}_T\subsetneq \mathcal{F}^{M}_T$. Moreover, $\mathbb{P}\left(\mathcal{T}^{(i)}(\Omega)\right)=1$. \\
2. From (\ref{gammabm=gammam}), it follows that $\mathcal{F}^{\mathcal{B}^{(i)}(M)}_T= \mathcal{F}^{M}_T$. Also, $\mathbb{P}\left(\mathcal{B}^{(i)}(\Omega)\right)=0$.  \end{rem}
\subsection{Connection between the cases $i=1$ and $i=2$}
Let $(Y_t)_{t\in[0,T]}$ be a Gaussian process with $Y_0=0$, a.s., and let \begin{equation}\mathcal{S}_t(Y)\ :=\  Y_T - Y_{T-t},\ t\in[0,T],\label{defS}\end{equation} be the time-reversion of the process $Y$.  Clearly, we have that $\mathcal{S}_0(Y)=0$, a.s., and \begin{equation}\mathcal{S}_T(Y)\ =\ Y_T,\ a.s.\label{S_T(M)=M_T}\end{equation}Furthermore, $\Gamma_T(\mathcal{S}(Y)) = \Gamma_T(Y)$ and $\mathcal{S}^2_t(Y)=Y_t$, a.s., $t\in[0,T]$.\\
\\
It is easy to show that $\mathcal{S}(M)$ is a continuous $\mathbb{F}_{T}^{\mathcal{S}(M)}$-martingale with \begin{equation}\langle \mathcal{S}(M)\rangle_{t}\ =\ \langle M\rangle_{T,T-t},\ t\in[0,T].\label{varianceofSM}\end{equation} Hence $\langle \mathcal{S}(M)\rangle_{\cdot}$ is strictly increasing on $[0,T]$. Note that $\mathcal{S}$ is a measure-preserving transformation if and only if $M$ has stationary increments, i.e. if and only if $M=W$.  By using (\ref{varianceofSM}), we obtain that 
 \begin{equation*}\left(\beta^{M,1}1_{[T-t,T)}\right)(\cdot)\ =\ \left(\beta^{\mathcal{S}(M),2}1_{[0,t)}\right)(T-\cdot),\ \lambda\text{-}a.e.,\ t\in(0,T].\label{b^1=bS2}\end{equation*}
It follows from this that \begin{equation}\mathcal{S}_t\left(\mathcal{T}^{(1)}(M)\right)\ =\ \mathcal{T}^{(2)}_t\left(\mathcal{S}(M)\right),\ a.s.,\ t\in[0,T].\label{ST=TS} \end{equation} 
Similarly, we have that  \begin{equation*}\left(\alpha^{M,1}1_{[T-t,T)}\right)(\cdot)\ =\ \left(\alpha^{\mathcal{S}(M),2}1_{[0,t)}\right)(T-\cdot),\ \lambda\text{-}a.e.,\ t\in(0,T].\label{a^1=aS2}\end{equation*}Therefore,
  \begin{equation}\mathcal{S}_t\left(\mathcal{B}^{(1)}(M)\right)\ =\ \mathcal{B}^{(2)}_t\left(\mathcal{S}(M)\right),\ a.s.,\ t\in[0,T].\label{SB=BS}\end{equation}
\section{Volterra Gaussian processes}
We  generalize the results of Section \ref{section2} to the case, when $X$ is a continuous Volterra Gaussian process with a non-degenerate Volterra kernel satisfying (\ref{kerinvertible}), (\ref{kernelsingular}) and that 
 \begin{equation}\mathbf{K}^{X}f\ =\  z_X(T,\cdot)1_{[0,t)} \text{ has a solution in } \check{\Lambda}_T(X)\text{ for all }t\in(0,T],\label{technicalassumption} \end{equation} which we denote by $k^{\ast}(t,\cdot):=k^{\ast}_{X,T}(t,\cdot)$, $t\in(0,T]$.
We proceed similarly as in Section \ref{section2}. 
\subsection{Measure-preserving and bridge transformations}
Let  $$M_t\  :=\ M_t(X,T)\ :=\ \text{E}_{\mathbb{P}}\left(X_T\big|\mathcal{F}_t^X\right), \ t\in[0,T],$$  denote the \textit{prediction martingale of $X_T$ with respect to $\mathbb{F}^X_T$}. By definition,
 \begin{equation}M_T\ =\ X_T,\ a.s. \label{X=M} \end{equation}
Clearly, by (\ref{filtrationX=filtrationW}) we have that \begin{equation}M_t\ =\ \int^t_0 z_X(T,s)dW_s,\ a.s.,\ t\in[0,T].\label{PM=intW}\end{equation}
From (\ref{PM=intW}) and (\ref{kerinvertible}), it follows that $M$ is a continuous Gaussian martingale with $M_0=0$, a.s., and the quadratic variation function $\langle M\rangle_{\cdot}=\int^{\cdot}_0 z^{2}_X(T,s)ds$ is strictly increasing on $[0,T]$. In particular, all results from Section \ref{section2} hold true for $M$. We have that $\Gamma_t(M)=\Gamma_t(W)$, $t\in[0,T]$, and hence from (\ref{assumption0}), it follows that \begin{equation}\Gamma_t(M)\ =\ \Gamma_t(X),\  t\in[0,T].\label{chaosM=chaosX}\end{equation}
More precisely, we
 have that\begin{equation}X_t\ =\ \int^t_0 k(t,s)dM_s,\ a.s.,\ t\in[0,T],\label{X=intkM}\end{equation} where the Volterra kernel is given by  \begin{equation}k(t,s)\ :=\ k_{X,T}(t,s)\ :=\ \frac{z_X(t,s)}{z_X(T,s)},\ s,t\in[0,T].\label{k=ztfraczT}\end{equation}  
 Let $\kappa:=\kappa^{X,T}$ be the isometric isomorphism defined by
 \begin{eqnarray*}\kappa\ :\  \Lambda_T(X)    &\to     & L^2_T(M)\\
                         1_{[0,t)}\ &\mapsto & k(t,\cdot),\ t\in(0,T].\end{eqnarray*} 
Then \begin{equation}\int^T_0 f(s)dX_s\ =\ \int^T_0 (\kappa f)(s)dM_s,\ a.s.,\ f\in\Lambda_T(X).\label{intfx=}\end{equation}
Also, from (\ref{PM=intW}) and (\ref{technicalassumption}), it follows that \begin{eqnarray}M_t&=& \int^T_0 k^{\ast}(t,s)dX_s,\ a.s.,\ t\in[0,T].\label{M=intkastX}\end{eqnarray}
Hence $$\kappa^{-1}1_{[0,t)}\ =\ k^{\ast}(t,\cdot),\ t\in(0,T].$$ From (\ref{chaosM=chaosX}), it follows that
 $k^{\ast}$ is a Volterra kernel.
Clearly, by using (\ref{X=M}), we have that \begin{equation}R^X(T,T)\ =\ \langle M\rangle_T\label{R^X=<M>}\end{equation} 
and \begin{equation}k^{\ast}(T,\cdot)\ \equiv\ 1.\ \label{kastT=1}\end{equation}From (\ref{k=ztfraczT}) and (\ref{kastT=1}), it follows that \begin{equation}\kappa(c)\ =\ \kappa^{-1}(c)\ =\ c,\ c\in\mathcal{C}_T.\label{kappa^{-1}c=c}\end{equation}
Moreover, (\ref{intfx=}) and (\ref{X=M}) imply that  \begin{eqnarray}\left(f,1_{[0,T)}\right)_{\Lambda_T(X)}
&=&\text{Cov}_{\mathbb{P}}\left(\int^T_0 f(s)dX_s, X_T\right)\nonumber\\
&=&\text{Cov}_{\mathbb{P}}\left(\int^T_0 (\kappa f)(s)dM_s, M_T\right)\nonumber\\
&=& \int^T_0 (\kappa f)(s)d\langle M\rangle_s,\ f\in\Lambda_T(X).\label{ettaXM}\end{eqnarray}
By combing (\ref{R^X=<M>}) and (\ref{ettaXM}), and then using (\ref{kappa^{-1}c=c}), we obtain that   \begin{equation}\eta^X f\ =\ \left(\kappa^{-1} \eta^M\kappa  \right)f,\ f\in\Lambda_T(X).\label{2eta}\end{equation} 
Also, by combining (\ref{pabr}), (\ref{2eta}), (\ref{intfx=}) and (\ref{notinvx}), we have that \begin{equation}\widehat{X}^T_t\ =\ \int^T_0 \overline{k(t,\cdot)}(s)d\widehat{M}^T_s,\ a.s., \ t\in[0,T].\label{antX=intantM}\end{equation}
\begin{exa}If $X=B^H$, then (\ref{technicalassumption}) is satisfied. We have that $$k^{\ast}(t,s)\ =\ 1_{[0,t)}(s) + \left(\frac{\sin\bigl(\pi\bigl(H-\frac{1}{2}\bigr)\bigr)}{\pi}s^{\frac{1}{2}-H}(t-s)^{\frac{1}{2}-H}\int^T_t \frac{u^{H-\frac{1}{2}}(u-t)^{{H-\frac{1}{2}}}}{u-s}du \right)1_{[0,t)}(s),$$ $s,t\in[0,T]$.
See \cite{Pi} for a full proof.\end{exa}
The following result generalizes Lemma \ref{pec}: 
\begin{lem}\label{genpec}
1. For $i\in\{1,2\}$ and $f\in\Lambda_T(X)$, let     \begin{equation}\alpha^{X,i}f\ :=\  \left(\kappa^{-1} \alpha^{M,i} \kappa\right) f\label{axi=}\end{equation} and  \begin{equation}\beta^{X,i}f\ :=\ \left(\kappa^{-1}\beta^{M,i}\kappa\right) f.\label{bxi=}\end{equation} 
Then $\alpha^{X,i}:\Lambda_{T,0}(X)\to\Lambda_T(X)$ and $\beta^{X,i}:\Lambda_T(X) \to\Lambda_{T,0}(X)$ are isometric isomorphisms with $\left(\alpha^{X,i}\right)^{-1}=\beta^{X,i}$. \\
2. For $i\in\{1,2\}$, we have that  \begin{equation}\alpha^{X,i}\eta^X f\ =\ \alpha^{X,i}f,\ f\in\Lambda_T(X),\label{genae=a}\end{equation}
and\begin{equation}\eta^X\beta^{X,i} f\ =\ \beta^{X,i}f,\ f\in\Lambda_T(X).\label{geneb=b}\end{equation}Also, \begin{equation}\beta^{X,i}\alpha^{X,i} f\ =\ \eta^X f,\ f\in\Lambda_T(X).\label{genba=e}\end{equation}\end{lem} 
\begin{proof}
1. Follows straightforward from part 2 of Lemma \ref{pec}.\\
2. Follows from part 3 of Lemma \ref{pec} by using (\ref{2eta}).
 \end{proof}
\begin{rem}\label{genoprem}1. Let $t\in(0,T]$ and $s\in[t,T)$. By combining part 1 of Remark \ref{oprem} and (\ref{kappa^{-1}c=c}), it follows that $\left(\alpha^{X,2}1_{[0,t)}\right)(s)=0$ and $\left(\beta^{X,2}1_{[0,t)}\right)(s)=c_M(k(t,\cdot),t)$.\\
2. From part 2 of Remark \ref{oprem} and (\ref{kappa^{-1}c=c}), we obtain that $\alpha^{X,i}c\equiv0$, $c\in\mathcal{C}_T$, $i\in\{1,2\}$.\end{rem}
The generalization of Theorem \ref{tablem} is straightforward: 
\begin{thm}\label{gentablem}Let $i\in\{1,2\}$. The transformation defined by \begin{equation}\mathcal{T}^{(i)}_t(X) \ :=\  \int^T_0 \left(\beta^{X,i}1_{[0,t)}\right)(s)dX_s,\ t\in[0,T],\label{mptx}\end{equation} is measure-preserving, i.e. $\mathcal{T}^{(i)}(X)$ is a centered Gaussian process with $\mathcal{T}^{(i)}_0(X)=0$, a.s., and $R^{\mathcal{T}^{(i)}(X)}=R^X$. Furthermore, the process 
\begin{equation}\mathcal{B}^{(i)}_t(X)\
:=\ \int^T_0 \left(\alpha^{X,i}1_{[0,t)}\right)(s)dX_s,\ t\in[0,T],\label{frtrans}\end{equation}
is a bridge of $X$.\end{thm}
\begin{proof}By construction, $\mathcal{T}^{(i)}(X)$ and $\mathcal{B}^{(i)}(X)$ are centered Gaussian processes.
By using part 1 of Lemma \ref{genpec}, we obtain that 
$$\text{Cov}_{\mathbb{P}}\left(\mathcal{T}^{(i)}_s(X),\mathcal{T}^{(i)}_t(X)\right)\ =\ \text{Cov}_{\mathbb{P}}(X_s,X_t),\ s,t\in[0,T].$$
Moreover, by combining (\ref{genae=a}), part 1 of Lemma \ref{genpec} and (\ref{pabr}), we obtain that 
$$\text{Cov}_{\mathbb{P}}\left(\mathcal{B}^{(i)}_s(X),\mathcal{B}^{(i)}_t(X)\right)\ =\ \text{Cov}_{\mathbb{P}}\left(X^T_s,X^T_t\right),\ s,t\in[0,T].\qedhere$$ \end{proof}
\begin{rem}Lemma \ref{genpec} and Theorem \ref{gentablem} hold true also without assumptions (\ref{kernelsingular}) and (\ref{technicalassumption}).\end{rem}
Next, we want to explicitly evaluate the functions $\alpha^{X,i}1_{[0,t)}$ and $\beta^{X,i}1_{[0,t)}$ in definitions (\ref{frtrans}) and (\ref{mptx}), respectively. For this purpose, we recall the notion of the Bochner integral in a Banach space: Let $(Q,\mathcal{Q},\mu)$ be a measure space with a finite signed measure and let $B$ be a Banach space with norm $|\cdot|_{B}$. The Bochner integral in $B$ of an indicator function $f(q):=g\cdot 1_{A}(q)$, $q\in Q$, where $A\in\mathcal{Q}$ and $g\in B$, is given by $\int_Q f(q) d\mu(q):=g\cdot\mu(A)$. A measurable map $f:(Q,\mathcal{Q})\to B$, $q\mapsto f(q)$, is Bochner integrable if $\int_{Q}|f(q)|_{B}d\mu(q)<\infty$. The space of Bochner integrable functions is a Banach space. The Bochner integral in $B$ of a Bochner integrable function is defined as the unique continuous, linear extension of $\int_{Q}\cdot\ d\mu$ from the set of indicator functions to $B$. By construction, the Bochner integral commutes with bounded linear maps: if $f:Q \to B$ is
Bochner integrable, $B'$ is a Banach space and $A:B\to B'$ is a bounded linear map, then $Af:Q\to B'$ is Bochner integrable, and \begin{equation*}A\int_Q f(q)d\mu(q)\ =\ \int_{Q}Af(q)d\mu(q).\label{commutes}\end{equation*} If the elements of $B$ are functions from $[0,T]$ to $\mathbb{R}$, then we denote $\int_Q f(q,s)d\mu(q):=\left(\int_{Q}f(q)d\mu(q)\right)(s)$, $s\in[0,T]$. For 
 details on  Bochner integration, see \cite{Mi}. \\
\\
We can state and prove the following lemma:
\begin{lem}\label{kappa^{-1}=}Let $z\in(0,T]$ and $g\in L^2_T(M)$.\\
1. If $g$ is absolutely continuous on $[a,z]$ for every $a>0$, then
\begin{equation*}\left(\kappa^{-1}g1_{[0,z)}\right)(s)\ =\ k^{\ast}(z,s)g(z)-\int^z_0 k^{\ast}(u,s)dg(u),\ \lambda\text{-}a.e.\ s\in(0,T),\label{rhoinv1}\end{equation*} where the integral is a Bochner integral in $\Lambda_T(X)$.\\
2. If $g$ is absolutely continuous on $[0,b]$ for every $b<z$, then
 \begin{equation*}\left(\kappa^{-1}g1_{[0,z)}\right)(s)\ =\ \int^z_0 \left(k^{\ast}(z,s)-k^{\ast}(u,s)\right)dg(u)+k^{\ast}(z,s)g(0),\  \lambda\text{-}a.e.\ s\in(0,T),\label{rhoinv3}\end{equation*} 
where the integral is a Bochner integral in $\Lambda_T(X)$.\end{lem}
\begin{proof}
1. We have that 
  \begin{equation*}\int^z_0 1_{[0,u)}(s)dg(u)\ =\ (g(z)-g(s))1_{[0,z)}(s),\ \lambda\text{-}a.e.\ s\in(0,T),\label{kappainvaz}\end{equation*} where the integral is a Bochner integral in $L^2_T(M)$. 
By combining this, the fact that $\kappa^{-1}:L^2_T(M)\to\Lambda_T(X)$ is a bounded linear map and Remark \ref{normf0ifff0}, we obtain that
 \begin{eqnarray*}\left(\kappa^{-1}g1_{[0,z)}\right)(s) &=&\left(\kappa^{-1}g(z) 1_{[0,z)}\right)(s)\ -\ \left(\kappa^{-1}(g(z)-g)1_{[0,z)}\right)(s)\nonumber\\&=& g(z)k^{\ast}(z,s)\ -\ \left(\kappa^{-1} \int^z_0 1_{[0,u)}dg(u)\right)(s)\nonumber\\
&=& g(z)k^{\ast}(z,s) \ -\ \left(\int^z_0 \kappa^{-1} 1_{[0,u)}dg(u)\right)(s)\nonumber\\
&=& g(z)k^{\ast}(z,s)\ -\ \int^z_0 k^{\ast}(u,s)dg(u),\ \lambda\text{-}a.e.\ s\in(0,T).\end{eqnarray*}
2. Similarly, we have that \begin{equation*}\int^z_0 \left(1_{[0,z)}(s)-1_{[0,u)}(s)\right)dg(u)\ =\ \left(g(s)\ -\ g(0)\right)1_{[0,z)}(s), \ \lambda\text{-}a.e.\ s\in(0,T),\label{kappainv0b}\end{equation*} where the integral is a Bochner integral in $L^2_T(M)$. Hence, \begin{eqnarray*}\left(\kappa^{-1}g1_{[0,z)}\right)(s)&=&\left(\kappa^{-1}\left(g-g(0)\right) 1_{[0,z)}\right)(s)\ +\ \left(\kappa^{-1}g(0)1_{[0,z)}\right)(s)\nonumber\\
&=&\left(\kappa^{-1}\int^z_0 \left(1_{[0,z)}-1_{[0,u)}\right)dg(u)\right)(s)\ +\ g(0)k^{\ast}(z,s)\nonumber\\
&=& \left(\int^z_0 \kappa^{-1}(1_{[0,z)}- 1_{[0,u)})dg(u)\right)(s)\ +\ g(0)k^{\ast}(z,s) \nonumber\\
&=& \int^z_0 \left(k^{\ast}(z,s)-k^{\ast}(u,s)\right)dg(u) + g(0)k^{\ast}(z,s),\ \lambda\text{-}a.e.\ s\in(0,T).\qedhere\end{eqnarray*}\end{proof}
Let $t\in(0,T]$. The function $\mathcal{H}^{M,1}k(t,\cdot)$ is absolutely continuous with respect to $\langle M\rangle_{\cdot}$, or equivalently, with respect to $\lambda(\cdot)$, on $[a,T]$ for every $a>0$. By combining part 1 of Lemma  \ref{kappa^{-1}=} with $z=T$ and (\ref{kastT=1}), we obtain that 
\begin{eqnarray}\left(\alpha^{X,1}1_{[0,t)}\right)(s)&=& 1_{[0,t)}(s)-\frac{1}{\langle M\rangle_T}\int^t_0 k(t,u)d\langle M\rangle_u\ +\
  \int^t_0 k^{\ast}(u,s) \frac{k(t,u)}{\langle M\rangle_u}d\langle M\rangle_u
\ \nonumber\\
&&-\int^T_0 \frac{k^{\ast}(u,s)}{\langle M\rangle^2_u}\int^u_0 k(t,v)d\langle M\rangle_v d\langle M\rangle_u,\ \lambda\text{-}a.e.\ s\in(0,T).\label{ax1=}\end{eqnarray}
Let $t\in(0,T)$. The function $\mathcal{H}^{M,2}k(t,\cdot)$ is absolutely continuous on $[0,T]$. Moreover,  $\left(\mathcal{H}^{M,2}k(t,\cdot)\right)(s)=0$, $s\geq t$. By using part 1 of Lemma \ref{kappa^{-1}=} with $z=t$, we obtain that
 \begin{eqnarray}\left(\alpha^{X,2}1_{[0,t)}\right)(s) &=& 
1_{[0,t)}(s)- \int^t_0 k^{\ast}(u,s)     \frac{k(t,u)}{\langle M\rangle_{T,u}}d\langle M\rangle_u\qquad \qquad\qquad\qquad\qquad\qquad\nonumber\\
& &+ \int^t_0 \frac{k^{\ast}(u,s)}{\langle M\rangle^2_{T,u}}\int^t_u k(t,v)d\langle M\rangle_v d\langle M\rangle_u,\ \lambda\text{-}a.e.\ s\in(0,T).\label{ax2=}\end{eqnarray}
Let $t\in(0,T]$. The function $\mathcal{H}^{M,1,\ast}k(t,\cdot)$ is absolutely continuous on $[a,T]$ for every $a>0$, and $\left(\mathcal{H}^{M,1,\ast}k(t,\cdot)\right)(s)=0$, $s\geq t$. By using part 1 of Lemma \ref{kappa^{-1}=} with $z=t$, we obtain that\begin{equation}\left(\beta^{X,1}1_{[0,t)}\right)(s)\ =\ 
1_{[0,t)}(s)- \int^t_0 k^{\ast}(u,s)\frac{k(t,u)}{\langle M\rangle_u} d\langle M\rangle_u,\ \lambda\text{-}a.e.\ s\in(0,T).\label{bx1=}\end{equation}
Let $t\in(0,T]$. The function $\mathcal{H}^{M,2,\ast}k(t,\cdot)$ is absolutely continuous on $[0,b]$ for every $b<T$. By using part 2 of Lemma \ref{kappa^{-1}=} with $z=T$ and (\ref{kastT=1}), we have that \begin{equation}\left(\beta^{X,2}1_{[0,t)}\right)(s)\ =\ 1_{[0,t)}(s)-
\int^t_0\left(1-k^{\ast}(u,s)\right)\frac{k(t,u)}{\langle M\rangle_{T,u}}d\langle M\rangle_u,\ \lambda\text{-}a.e.\ s\in(0,T).\label{bx2=}\end{equation}
\subsection{Pathwise relations for the transformed processes}
By combining (\ref{mptx}), (\ref{bxi=}), (\ref{intfx=}) and (\ref{pwc2}), we obtain that
\begin{eqnarray}\mathcal{T}^{(i)}_t(X)&=&
\int^t_0 k(t,s)d\mathcal{T}^{(i)}_s(M),\ a.s.,\ t\in[0,T],\ i\in\{1,2\}
,\label{TX=intTM}\end{eqnarray}i.e. $\mathcal{T}^{(i)}(M)$ is the prediction martingale of $\mathcal{T}^{(i)}_T(X)$ with respect to $\mathbb{F}_T^{\mathcal{T}^{(i)}(X)}$. Clearly, we can evaluate $\mathcal{T}^{(i)}(X)$ explicitly by using (\ref{bx1=}) and (\ref{bx2=}). However, in order to obtain expressions corresponding to (\ref{genyor}) and (\ref{T^2(M)=}), 
we combine (\ref{mptx}), (\ref{intfx=}), the stochastic Fubini theorem and (\ref{M=intkastX}), and obtain that 
\begin{equation*}\mathcal{T}^{(1)}_t(X)\ =\ X_t -\ \int^t_0 \left(\int^u_0  k^{\ast}(u,s)dX_s \right)\frac{k(t,u)}{\langle M\rangle_u}d\langle M\rangle_u,\ a.s.,\ t\in[0,T].\end{equation*}
In the same way, we have that \begin{equation}\mathcal{T}^{(2)}_t(X)\ =\ X_t-\int^t_0\left(X_T-\int^u_0 k^{\ast}(u,s)dX_s\right)\frac{k(t,u)}{\langle M\rangle_{T,u}}d\langle M\rangle_u,\ a.s.,\ t\in[0,T].\label{evT^2(X)=intX}\end{equation}
On the one hand, by combining (\ref{mptx}), (\ref{geneb=b}) and (\ref{notinvx}), we have that \begin{eqnarray}\mathcal{T}^{(i)}_t(X)
&=&\int^T_0 \left(\overline{\beta^{X,i} 1_{[0,t)}}\right)(s)d\widehat{X}^T_s,\ a.s., \ t\in[0,T],\ i\in\{1,2\}.\label{genT=intAB}\end{eqnarray}
On the other hand, by combining (\ref{pabr}), (\ref{genba=e}) and (\ref{mptx}), we obtain the  reciprocal \begin{eqnarray}\widehat{X}^T_t 
&=& \int^T_0 \left(\alpha^{X,i}1_{[0,t)}\right)(s)d\mathcal{T}^{(i)}_s(X),\ a.s.,\ t\in[0,T],\ i\in\{1,2\}.\label{genAB=intT}\end{eqnarray}
From (\ref{genT=intAB}) and (\ref{genAB=intT}), we have that \begin{equation}\Gamma_T\!\left(\mathcal{T}^{(i)}(X)\right)\ =\ \Gamma_T\!\left(\widehat{X}^T\right),\ i\in\{1,2\}.\label{gammatx=gammaantx}\end{equation}
For $i=2$, it follows from part 1 of Remark \ref{genoprem} that $$\Gamma_t\!\left(\mathcal{T}^{(2)}(X)\right)\ =\ \Gamma_t\!\left(\widehat{X}^T\right),\ t\in[0,T],$$ i.e. $$\mathbb{F}^{\mathcal{T}^{(2)}(X)}_{T}\ =\ \mathbb{F}^{\widehat{X}^T}_{T}.$$
By combining (\ref{evT^2(X)=intX}), (\ref{intX=intantX+R}) and (\ref{ettaXM}), we obtain that
 \begin{equation*}\mathcal{T}^{(2)}_t(X)\ =\ \widehat{X}^T_t\ +\   
\int^t_0\int^T_0 \overline{k^{\ast}(u,\cdot)}(s)d\widehat{X}^T_s\frac{k(t,u)}{\langle M\rangle_{T,u}}d\langle M\rangle_u,\ a.s.,\ t\in[0,T].\label{ANTB=intT2}\end{equation*} Note that, in contrast to $\widehat{M}^T$, the process $\widehat{X}^T$ is \textit{not} the solution of a linear stochastic differential equation (unless $X=M$).  
By combining (\ref{frtrans}), (\ref{axi=}), (\ref{intfx=}) and (\ref{pwc1}), we obtain that  \begin{eqnarray}\mathcal{B}^{(i)}_t(X)&=&\int^T_0 \overline{k(t,\cdot)}(s)d\mathcal{B}^{(i)}_s(M),\ a.s.,\ t\in[0,T],\ i\in\{1,2\}.\label{BX=intBM}\end{eqnarray}
By using (\ref{ax1=}), we obtain that  
\begin{eqnarray*}\mathcal{B}^{(1)}_t(X)& =&  X_t-\frac{X_T}{\langle M\rangle_T}\int^t_0 k(t,u)d\langle M\rangle_u\ +\ \int^t_0\left( \int^t_0 k^{\ast}(u,s)\frac{k(t,u)}{\langle M\rangle_u}d\langle M\rangle_u\right)dX_s\\
&&-\int^T_0\left(\int^T_0  \frac{k^{\ast}(u,s)}{\langle M\rangle_u^2}\int^u_0 k(t,v)d\langle M\rangle_vd\langle M\rangle_u\right)dX_s,\ a.s.,\ t\in(0,T].\end{eqnarray*}
Similarly, by using (\ref{ax2=}), we have that 
 \begin{eqnarray*}\mathcal{B}^{(2)}_t(X)&=& X_t -
\int^t_0 \left(\int^t_0 k^{\ast}(u,s)     \frac{k(t,u)}{\langle M\rangle_{T,u}}d\langle M\rangle_u\right) dX_s\nonumber\\
& &+ \int^t_0 \left(\int^t_0\frac{k^{\ast}(u,s)}{\langle M\rangle^2_{T,u}}\int^t_u k(t,v)d\langle M\rangle_v d\langle M\rangle_u\right)dX_s,\ a.s.,\ t\in[0,T).
\end{eqnarray*}
By combining part 1 of Lemma \ref{genpec} and (\ref{frtrans}),
 we have that \begin{equation}X_t\ =\ \int^T_0 \left(\overline{\beta^{X,i}1_{[0,t)}}\right)(s)d\mathcal{B}^{(i)}_s(X),\ a.s.,\ t\in[0,T], \ i\in\{1,2\}.\label{X=intBiM}\end{equation}
By comparing identities (\ref{frtrans}) and (\ref{X=intBiM}) with identities (\ref{genAB=intT}) and (\ref{genT=intAB}), we conclude that the bridge $\mathcal{B}^{(i)}(X)$ is related to the process $X$ in the same way, as 
the bridge $\widehat{X}^T$ is related to the process $\mathcal{T}^{(i)}(X)$. Therefore, \begin{equation}\Gamma_T\left(\mathcal{B}^{(i)}(X)\right)\ =\ \Gamma_T(X),\ i\in\{1,2\}.\label{gammabx=gammax}\end{equation}
Furthermore, $$\Gamma_t\left(\mathcal{B}^{(2)}(X)\right)\ =\ \Gamma_t(X),\ t\in[0,T],$$
i.e. $$\mathbb{F}^{\mathcal{B}^{(2)}(X)}_T\ =\ \mathbb{F}^X_T.$$\begin{rem}$\mathcal{B}^{(2)}(X)$ is the dynamic bridge of $X$, and has been introduced in \cite{Ga}.\end{rem}
From (\ref{genAB=intT}) and (\ref{X=intBiM}), we obtain the following: 
\begin{lem}\label{lemtt}Let $i\in\{1,2\}$. Then
 \begin{equation*}  \widehat{X}_t^{T} =\ \mathcal{B}^{(i)}_t\left(\mathcal{T}^{(i)}(X)\right),   \ a.s., \ t\in[0,T].\end{equation*}
Furthermore,  $$X_t\ =\ \mathcal{T}^{(i)}_t\left(\mathcal{B}^{(i)}(X)\right),\ a.s.,\ t\in[0,T].$$\end{lem}
\begin{rem}1. From (\ref{gammatx=gammaantx}), we have that $\mathcal{F}^{\mathcal{T}^{(i)}(X)}_T\subsetneq \mathcal{F}^{X}_T$. Moreover, $\mathbb{P}\left(\mathcal{T}^{(i)}(\Omega)\right)=1$. \\
2. From (\ref{gammabx=gammax}), it follows that $\mathcal{F}^{\mathcal{B}^{(i)}(X)}_T= \mathcal{F}^{X}_T$. Also, $\mathbb{P}\left(\mathcal{B}^{(i)}(\Omega)\right)=0$. \end{rem}
\subsection{Connection between the cases $i=1$ and $i=2$}
Consider (\ref{X=intkM}) and let \begin{eqnarray*}X^{\mathcal{S}}_t& :=& \int^t_0 k(t,s)d\mathcal{S}_{s}(M)\\
&=&\int^t_0 \left(\frac{z_X(t,s)}{z_X(T,s)}z_X(T,T-s)\right)d\mathcal{S}_s(W),\ a.s.,\ t\in[0,T],\end{eqnarray*} where $\mathcal{S}$ is defined  as in (\ref{defS}). Note that $X^{\mathcal{S}}$ is well-defined by part 1 of Remark \ref{canonicalrepresentation}. Clearly, $X^{\mathcal{S}}$ is a continuous Volterra Gaussian process with a non-degenerate Volterra kernel, and $\mathcal{S}(M)$ is the prediction martingale of $X^{\mathcal{S}}_T$ with respect to $\mathbb{F}^{X^{\mathcal{S}}}_T=\mathbb{F}_T^{\mathcal{S}(M)}$. 
By combining (\ref{TX=intTM}), (\ref{ST=TS}) and again (\ref{TX=intTM}), we obtain that \begin{equation*}\left(\mathcal{T}^{(1)}(X)\right)^{\mathcal{S}}_t\ =\ \mathcal{T}^{(2)}_t\left(X^{\mathcal{S}}\right),\ a.s.,\ t\in[0,T]. \label{TnonVolterra}\end{equation*}
Clearly,  \begin{equation}R^{\mathcal{S}(X)}(T,T-t)\ =\ R^X(T,T)-R^X(T,t),\ t\in[0,T]. \label{SANTB=ANTBS}\end{equation}By using (\ref{SANTB=ANTBS}) and (\ref{S_T(M)=M_T}), we obtain that    
 \begin{equation} 
\mathcal{S}_t\left(\widehat{X}^T\right)\ =\ \widehat{\mathcal{S}(X)}^T_t,\ a.s.,\ t\in[0,T].\label{SANTBX=ANTBSX}\end{equation}
Consider (\ref{antX=intantM}) and let
\begin{equation}\left(\widehat{X}^T\right)^{\mathcal{S}}_t\ :=\ \int^T_0 \overline{k(t,\cdot)}(s)d\mathcal{S}_s\left(\widehat{M}^T\right),\ t\in[0,T].\label{XS=intkSANTM}\end{equation} Due to (\ref{SANTBX=ANTBSX}), the right-hand side of (\ref{XS=intkSANTM}) is a Wiener integral with respect to $\widehat{\mathcal{S}(M)}^T$. The process $\left(\widehat{X}^T\right)^{\mathcal{S}}$ is well-defined  due to part 1 of Remark \ref{canonicalrepresentation}. In particular, from (\ref{antX=intantM}) it follows that $$\left(\widehat{X}^T\right)^{\mathcal{S}}_t\ =\ \widehat{X^\mathcal{S}_t}^T,\ a.s.,\ t\in[0,T].$$  
By combining (\ref{BX=intBM}), (\ref{SB=BS}) and again (\ref{BX=intBM}), we have that \begin{equation*}\left(\mathcal{B}^{(1)}(X)\right)_t^{\mathcal{S}} \ =\ \mathcal{B}^{(2)}_t\left(X^{\mathcal{S}}\right),\ a.s.,\ t\in[0,T].\label{BnonVolterra}\end{equation*}
\begin{rem} If $X$ has stationary increments, then $\mathcal{S}(X)\stackrel{d}{=}X$, and hence $\mathcal{S}(X)$ is a Volterra Gaussian process. In general, however, it is not clear whether $\mathcal{S}(X)$ is a Volterra Gaussian process.
 \end{rem}
\section{A Fourier-Laguerre series expansion}
Let $(M_t)_{t\in[0,\infty)}$ be a continuous Gaussian martingale with $M_0=0$, a.s., such that $\langle M\rangle_{\cdot}$ is strictly increasing and \begin{equation}\lim_{t\to\infty}\langle M\rangle_{t}\ =\ \infty.\label{condtimechange}\end{equation}
As before, we assume that the underlying probability space is the coordinate space of $M$, i.e. we assume that $\Omega=\{\omega:[0,\infty)\to\mathbb{R}\,|\,\omega\text{ is continuous}\}$, $\mathcal{F}=\mathcal{F}^M_{\infty}:=\sigma(M_t\,|\,t\in[0,\infty))$ and $\mathbb{P}$ is the probability measure with respect to which $M_t(\omega)=\omega(t)$, $\omega\in\Omega$, $t\in[0,\infty)$, is a Gaussian martingale with quadratic variation function $\langle M\rangle_{\cdot}$.
Let $$\Gamma_{\infty}(M)\ :=\ \overline{\text{span}\{M_t\,|\,t\in[0,\infty)\}}$$ and $$\Gamma_{[T,\infty)}(M)\ :=\ \overline{\text{span}\{M_t-M_s\,|\, s,t\in[T,\infty)\}}$$ denote the first Wiener chaoses of $M$ over $[0,\infty)$ and $[T,\infty)$, respectively. By (\ref{condtimechange}), there exists a standard Brownian motion $(W_t)_{t\in[0,\infty)}$ such that $M_t=W_{\langle M\rangle_t}$, $a.s.$, $t\in[0,\infty)$ (see \cite{Ka}, Theorem 4.6, p. 174). From the strong law of large numbers for standard Brownian motion, it follows that\begin{eqnarray}\lim_{t\to\infty}\frac{M_t}{\langle M\rangle_t}\ =\ \lim_{t\to\infty}\frac{W_t}{t}\ =\  0,\ a.s.\label{SLLN}\end{eqnarray}
From (\ref{genyor}), we see that transformation $\mathcal{T}^{(1)}$ does not depend on $T$, hence we can write \begin{eqnarray*}\mathcal{T}_t(M)&:=& \mathcal{T}^{(1)}_t(M)\ =\ M_t- \int^t_0\frac{M_s}{\langle M\rangle_s}d\langle M\rangle_s,\ a.s., \ t\in[0,\infty).\end{eqnarray*}
By using the stochastic Fubini theorem, we have that
\begin{eqnarray*}\mathcal{T}_t(M)&=&   \int^{\infty}_0 \left(\beta_{\infty}^{M,1}1_{[0,t)}\right)(s)dM_s,\ a.s.,\ t\in[0,\infty),\end{eqnarray*}
where 
\begin{eqnarray*}\beta^{M,1}_{\infty}\ :\ L^2([0,\infty),d\langle M\rangle_{\cdot})&\to&  L^2([0,\infty),d\langle M\rangle_{\cdot})\\
f(\cdot)\quad\qquad&\mapsto & f(\cdot)-\int^{\infty}_{\cdot} \frac{f(u)}{\langle M\rangle_u}d\langle M\rangle_u.\end{eqnarray*} By combining the classical Fubini theorem and (\ref{condtimechange}), we obtain that $\beta^{M,1}_{\infty}$ is an isometric isomorphism with inverse 
 \begin{eqnarray*}\alpha^{M,1}_{\infty}\ :\ L^2([0,\infty),d\langle M\rangle_{\cdot})&\to&  L^2([0,\infty),d\langle M\rangle_{\cdot})\\
f(\cdot)\quad\qquad&\mapsto & f(\cdot)-\frac{1}{\langle M\rangle_{\cdot}}\int^{\cdot}_0 f(u)d\langle M\rangle_u.\end{eqnarray*}
Hence, the measure-preserving transformation $\mathcal{T}$ is an automorphism, i.e. it has a measurable inverse. The inverse is given by \begin{equation}\mathcal{T}^{-1}_t(M)\ =\ \int^{\infty}_0\left(\alpha^{M,1}_{\infty}1_{[0,t)}\right)(s)dM_s
\  =\ -\langle M\rangle_t \int^{\infty}_t \frac{1}{\langle M\rangle_s}dM_s,\ a.s.,\ t\in(0,\infty).
 \label{T^{-1}=intadM}\end{equation} In particular, \begin{equation}\Gamma_{\infty}\left(\mathcal{T}^{-1}(M)\right)\ =\ \Gamma_{\infty}(M).\label{chaoseq3}\end{equation}
\begin{rem}For the (on $T$ dependent) bridge of $M$ in (\ref{B1M=}), it holds by using (\ref{T^{-1}=intadM}) that $$\lim_{\begin{tiny}\begin{array}{c}T\to\infty\\
T>t\end{array}\end{tiny}}\mathcal{B}^{(1)}_t(M)\ =\ \mathcal{T}^{-1}_t(M),\ a.s.,\ t\in(0,\infty).$$\end{rem}
Let $T>0$ be fixed. $M$ has independent increments, hence \begin{equation}\Gamma_{\infty}(M)\ =\ \Gamma_T(M)\ \perp\ \Gamma_{[T,\infty)}(M).\label{chaoseq}\end{equation}
On the one hand, from (\ref{notsamex}) and (\ref{GT=GANTB}), it follows that \begin{equation}\Gamma_T(M)\ =\ \Gamma_T\left(\mathcal{T}(M)\right)\ \perp\ \text{span}\{M_T\}.\label{firstit}\end{equation}
On the other hand, by using (\ref{chaoseq3}) and (\ref{firstit}) with $\mathcal{T}^{-1}(M)$ instead of $M$,  we have that 
\begin{eqnarray}\Gamma_{[T,\infty)}(M)&= &
\Gamma_{[T,\infty)}\left(\mathcal{T}^{-1}(M)\right)\ \perp\  \text{span}\left\{\mathcal{T}^{-1}_T(M)\right\}.\label{secondit}\end{eqnarray} 
By combining (\ref{chaoseq}), (\ref{firstit}) and  (\ref{secondit}), and iterating this decomposing procedure, we obtain that
 $\left\{\mathcal{T}^n_T(M)\right\}_{n\in\mathbb{Z}}$ is an orthogonal system in $\Gamma_{\infty}(M)$. 
From the following two-sided Fourier-Laguerre series expansion, it follows that this system is complete:
\begin{thm}Let $Z:=\int^{\infty}_0 f(s)dM_s\in\Gamma_{\infty}(M)$. Then 
$$Z\ =\ L^2(\mathbb{P})\text{-}\sum_{n\in\mathbb{Z}}\left(\int^{\infty}_0 f(s)\mathcal{L}^{M,T}_n(s)d\langle M\rangle_s\right)\cdot \varepsilon^{M,T}_n,$$
where the sequence $\left\{\varepsilon_n^{M,T}\right\}_{n\in\mathbb{Z}}:=\left\{\frac{\mathcal{T}^n_T(M)}{\sqrt{\langle M\rangle_T}}\right\}_{n\in\mathbb{Z}}$ is  i.i.d. with $\varepsilon^{M,T}_0\sim \mathcal{N}(0,1)$, and 
\begin{equation*}\mathcal{L}^{M,T}_n(s)\ :=\begin{cases} \frac{1}{\sqrt{\langle M\rangle_T}} L_n\left(\ln\frac{\langle M\rangle_{T}}{\langle M\rangle_{s}}\right)1_{(0,T)}(s), \qquad & n\in\mathbb{N}_0\\
\\
 -\frac{\sqrt{\langle M\rangle_T}}{\langle M\rangle_s} L_{-n-1}\left(\ln\frac{\langle M\rangle_{s}}{\langle M\rangle_{T}}\right)1_{(T,\infty)}(s), & n\in -\mathbb{N},\end{cases}\end{equation*}  with 
 $L_n(x):=\sum_{k=0}^n \binom{n}{k}\frac{1}{k!}(-x)^k$ denoting the $n$-th Laguerre polynomial, $n\in\mathbb{N}_0$. \label{seriesexpansion}\end{thm}
\begin{proof}
First, we show that $\left\{\varepsilon^{M,T}_n\right\}_{n\in\mathbb{N}_0}$ is a Hilbert basis of $\Gamma_T(M)$. By iterating $\mathcal{T}$ and using the stochastic Fubini theorem, we have that 
\begin{eqnarray}\mathcal{T}^n_T(M)&=&  M_T+  \sum_{k=1}^{n}(-1)^k \binom{n}{k}\int^T_0 \int^{t_{k-1}}_0\ldots \int^{t_1}_0 M_s\frac{d\langle M\rangle_s}{\langle M\rangle_s}\frac{d\langle M\rangle_{t_1}}{\langle M\rangle_{t_1}}\dots \frac{d\langle M\rangle_{t_{k-1}}}{\langle M \rangle_{t_{k-1}}}\nonumber \\
&=& M_T +  \sum_{k=1}^{n}\frac{(-1)^k}{k!} \binom{n}{k}\int^T_0 \ln^k\left(\frac{\langle M\rangle_T}{\langle M\rangle_s}\right)\nonumber dM_s\\
&=& \int^T_0 L_n\left(\ln\frac{\langle M\rangle_T}{\langle M\rangle_s}\right)dM_s\nonumber\\
&=&\sqrt{\langle M\rangle_T}\int^T_0\mathcal{L}^{M,T}_n(s)dM_s,\ a.s.,\ n\in\mathbb{N}_0.\label{iterationTn}\end{eqnarray}
It is well-known that \begin{equation}\{L_n(x),x\in[0,\infty)\}_{n\in\mathbb{N}_0}\text{ is a Hilbert basis 
of }L^2([0,\infty),e^{-x}dx).\label{laguerrebasis}\end{equation} By substitution, we obtain that $\left\{\mathcal{L}^{M,T}_n(s),s\in(0,T)\right\}_{n\in\mathbb{N}_0}$  is a Hilbert basis of $L^2_T(M)$. The claim follows from (\ref{iterationTn}) and the Wiener isometry between $L^2_T(M)$ and $\Gamma_T(M)$.
Second, we show that $\left\{\varepsilon^{M,T}_n\right\}_{n\in-\mathbb{N}}$ is a Hilbert basis of $\Gamma_{[T,\infty)}(M)$. 
By using (\ref{T^{-1}=intadM}), partial integration and (\ref{SLLN}), we have that  \begin{equation}\mathcal{T}^{-1}_t(M)\ =\ M_t\ -\ \langle M\rangle_t \cdot \mathcal{J}^1_t(M),\ a.s., \ t\in(0,\infty),\label{bl4}\end{equation} where
$$\mathcal{J}^k_t(M):=\int^{\infty}_t \int^{\infty}_{t_{k-1}}\ldots \int^{\infty}_{t_1} \frac{M_s}{\langle M\rangle_s}\frac{d\langle M\rangle_s}{\langle M\rangle_s}\frac{d\langle M\rangle_{t_1}}{\langle M\rangle_{t_1}}\dots \frac{d\langle M\rangle_{t_{k-1}}}{\langle M \rangle_{t_{k-1}}},\ k\in\mathbb{N}.$$ 
By combining (\ref{bl4}) and (\ref{T^{-1}=intadM}), using the stochastic Fubini theorem and iterating, we obtain that
\begin{eqnarray}\mathcal{J}^k_T(M)&=&\int^{\infty}_{T}\int^{\infty}_{t_{k-1}}\ldots\int^{\infty}_{t_2}\left(\frac{M_{t_1}}{\langle M\rangle_{t_1}}+\int^{\infty}_{t_1}\frac{dM_s}{\langle M\rangle_s}\right)\frac{d\langle M\rangle_{t_1}}{\langle M\rangle_{t_1}}\cdots \frac{d\langle M\rangle_{t_{k-1}}}{\langle M \rangle_{t_{k-1}}}\nonumber\\
&=&\mathcal{J}^{k-1}_T(M)\ +\ \int^{\infty}_T \frac{1}{(k-1)!}\ln^{k-1}\left(\frac{\langle M\rangle_s}{\langle M\rangle_T}\right)\frac{dM_s}{\langle M\rangle_s}\nonumber\ =\ \ldots\\
&=& \mathcal{J}^1_T(M)\ +\ \sum_{j=1}^{k-1}\int^{\infty}_T \frac{1}{j!} \ln^{j}\left(\frac{\langle M\rangle_s}{\langle M\rangle_T}\right)\frac{dM_s}{\langle M\rangle_s},\ a.s.,\  k\geq 2.\label{bl3}\end{eqnarray}
From the identity $\binom{n}{k}=\binom{n-1}{k}+\binom{n-1}{k-1}$, $1< k\leq n$, it follows that  \begin{equation}\sum_{k=j+1}^{n}(-1)^{k+j-1}\binom{n}{k}\ =\ \binom{n-1}{j},\ 0\leq j<n.\label{bl2}\end{equation}
By using (\ref{bl4}) and iterating, then using (\ref{bl3}), again (\ref{bl4}) and (\ref{bl2}), we obtain that
\begin{eqnarray}\mathcal{T}^{-n}_T(M)&=& M_T+\sum_{k=1}^{n}(-1)^k \binom{n}{k}\langle M\rangle_T\cdot \mathcal{J}^k_T(M)\nonumber\\
&=&M_T+\sum_{k=1}^{n}(-1)^k \binom{n}{k}\langle M\rangle_T \left(\mathcal{J}^1_T(M)+ \sum_{j=1}^{k-1}\int^{\infty}_T \frac{1}{j!}\ln^{j}\left(\frac{\langle M\rangle_s}{\langle M\rangle_T}\right)\frac{dM_s}{\langle M\rangle_s}\right)\nonumber\end{eqnarray}
\begin{eqnarray}&=&\sum_{k=1}^{n}(-1)^k \binom{n}{k}\langle M\rangle_T  \left(\int^{\infty}_T\frac{1}{\langle M\rangle}_s dM_s\ +\  \sum_{j=1}^{k-1}\int^{\infty}_T \frac{1}{j!}\ln^{j}\left(\frac{\langle M\rangle_s}{\langle M\rangle_T}\right)\frac{dM_s}{\langle M\rangle_s}\right)\nonumber\\
&=& \langle M\rangle_T \int^{\infty}_T \left( -1 + \sum_{k=1}^{n}(-1)^k \binom{n}{k} \sum_{j=1}^{k-1} \frac{1}{j!}\ln^{j}\left(\frac{\langle M\rangle_s}{\langle M\rangle_T}\right)\right) \frac{dM_s}{\langle M\rangle_s}\nonumber\\
&=&
-\langle M\rangle_T \int^{\infty}_T \left( 1 +\sum_{j=1}^{n-1}\sum_{k=j+1}^n (-1)^{k-1+j} \binom{n}{k} \frac{1}{j!}\ln^j \left(\frac{\langle M\rangle_s}{\langle M\rangle_T}\right) (-1)^j  \right) \frac{dM_s}{\langle M\rangle_s}\nonumber\\
&=& - \langle M\rangle_T\int^{\infty}_{T} L_{n-1}\left(\ln\frac{\langle M\rangle_s}{\langle M\rangle_T}\right)\frac{dM_s}{\langle M\rangle_s}\nonumber\\
&=& \sqrt{\langle M\rangle_T}\int^{\infty}_T \mathcal{L}_{-n}^{M,T}(s)dM_s,\ a.s.,\ n\in\mathbb{N}.\label{T-1=iteration}\end{eqnarray}
From (\ref{laguerrebasis}), we obtain by substitution that $\left\{\mathcal{L}^{M,T}_{n}(s),s\in(T,\infty)\right\}_{n\in -\mathbb{N}}$ is a Hilbert basis of $L^2([T,\infty),d\langle M\rangle_{\cdot})$. By combining (\ref{T-1=iteration}) and the Wiener isometry between $L^2([T,\infty),d\langle M\rangle_{\cdot})$ and $\Gamma_{[T,\infty)}(M)$, we obtain the claim.
Third, by combining (\ref{chaoseq}) with these results, we obtain that  $\left\{\varepsilon^{M,T}_n\right\}_{n\in\mathbb{Z}}$ is a Hilbert basis of $\Gamma_{\infty}(M)$. A Fourier expansion yields \begin{eqnarray*}Z &=& L^2(\mathbb{P})\text{-}\sum_{n\in\mathbb{Z}} \text{Cov}_{\mathbb{P}}\left(Z, \varepsilon^{M,T}_n\right)\cdot \varepsilon^{M,T}_n\\
&=& L^2(\mathbb{P})\text{-}\sum_{n\in\mathbb{Z}}\left(\int^{\infty}_0 f(s)\mathcal{L}^{M,T}_n(s)d\langle M\rangle_s\right) \cdot\varepsilon_n^{M,T}.\qedhere\end{eqnarray*}\end{proof}
As a special case, we obtain the following one-sided Fourier-Laguerre series expansion. For $M=W$, it was shown in \cite{Je}.
\begin{cor}We have that $$M_t\ =\ L^2(\mathbb{P})\text{-}\sum_{n\in\mathbb{N}_0}\sqrt{\langle M\rangle_T}\cdot\zeta_n\left(\ln\frac{\langle M\rangle_T}{\langle M\rangle_t}\right)\cdot\varepsilon_n^{M,T},\ t\in(0,T],$$
where $\zeta_n(y):=\int^{\infty}_y L_n(x)e^{-x}dx$, $y\in[0,\infty).$
\end{cor}\begin{proof}Set $f:=1_{[0,t)}$ in Theorem \ref{seriesexpansion}.
\end{proof} 
\begin{rem}
Clearly, we have that \begin{equation*}\mathcal{F}^M_T\ =\ \bigvee_{n\in\mathbb{N}_0}\sigma\left(\mathcal{T}^n_T(M)\right)\label{F^M_T} \end{equation*} and \begin{equation*}\mathcal{F}^M_{[T,\infty)}\  :=\ \sigma(M_t-M_s\,|\,s,t\in[T,\infty))\ =\ \bigvee_{n\in-\mathbb{N}}\sigma\left(\mathcal{T}^n_T(M)\right).\label{G^M_T} \end{equation*}
Furthermore, \begin{equation*}\mathcal{F}\ =\ \bigvee_{n\in\mathbb{Z}}\sigma\left(\mathcal{T}^{n}_T(M)\right).\end{equation*}
Recall that an automorphism $\mathcal{T}$ on $(\Omega,\mathcal{F},\mathbb{P})$ is a \textit{Kolmogorov automorphism}, if there exists a $\sigma$-algebra $\mathcal{A}\subseteq\mathcal{F}$, such that  $\mathcal{T}^{-1}\mathcal{A}\subseteq \mathcal{A}$, $\vee_{n\in\mathbb{Z}}\mathcal{T}^n\mathcal{A} =\mathcal{F}$ and $\cap_{n\in\mathbb{N}_0} \mathcal{T}^{-n}\mathcal{A}=\{\Omega,\emptyset\}$.
It is straightforward to see that $\mathcal{T}$ and $\mathcal{T}^{-1}$ are Kolmogorov automorphisms with $\mathcal{A}=\mathcal{F}^M_{[T,\infty)}$ and 
$\mathcal{A}=\mathcal{F}^M_{T}$, respectively. 
Hence, $\mathcal{T}$ and $\mathcal{T}^{-1}$ are strongly mixing and hence ergodic (see \cite{Pet}, Proposition 5.11 and Proposition 5.9 on p. 63 and p. 62).\end{rem}
\textbf{Acknowledgements.} Thanks are due to my supervisor Esko Valkeila for motivating the writing of this article, for good questions and helpful comments. I thank Giovanni Peccati for a particularly motivating discussion. Also, I thank Ilkka Norros for helpful advices. I am indebted to the Finnish Graduate School in Stochastics (FGSS) for financial support.

\end{document}